\documentclass[12pt]{amsart}

\usepackage{latexsym}
\usepackage{amssymb}
\usepackage[cp850]{inputenc}
\usepackage{epsfig}
\usepackage{psfrag}
\usepackage{amsthm}
\usepackage{amscd}
\usepackage{amsmath}
\usepackage{amsfonts}
\usepackage{graphics}
\usepackage[all]{xy}

\newtheorem{theorem}{Theorem}[section] 
\newtheorem{lemma}[theorem]{Lemma}
\newtheorem{cor}[theorem]{Corollary} 
\newtheorem{prop}[theorem]{Proposition}

\newtheorem*{theorem*}{Theorem} 
\newtheorem*{corollary*}{Corollary}

\theoremstyle{definition}
\newtheorem{example}{Example}[section]

\newtheorem{definition}[example]{Definition}
\newtheorem*{remark*}{Remark}

\newtheorem*{definition*}{Definition}
\newtheorem*{example*}{Example}

\theoremstyle{theorem}
\let\ssection=\section
\renewcommand{\section}{\setcounter{equation}{0}\ssection}

\newtheorem*{namedtheorem}{\theoremname}
\newtheorem*{maintheoremenv}{Theorem 1.1}

\newtheorem*{boundeddiameterlemma}{Bounded Diameter Lemma}
\theoremstyle{remark}

\newcommand{\BC}{\mathbb C}			
\newcommand{\BR}{\mathbb R}			
			
			\newcommand{\BZ}{\mathbb Z}

\newcommand{\CG}{\mathcal G}		%\newcommand{\CH}{\mathcal H}

		\newcommand{\CN}{\mathcal N}

\DeclareMathOperator{\PSL}{PSL}		%	Spezielle lineare Gruppe
\DeclareMathOperator{\vol}{vol}		%	Volumen

\DeclareMathOperator{\waist}{waist}

\DeclareMathOperator{\inj}{inj}
\DeclareMathOperator{\axis}{axis}

\DeclareMathOperator{\diam}{diam}
\DeclareMathOperator{\rank}{rank}

\DeclareMathOperator{\length}{length}

\DeclareMathOperator{\dist}{d}

\DeclareMathOperator{\CH}{CH}
\DeclareMathOperator{\TCH}{TCH}

\DeclareMathOperator{\CC}{CC}

\newcommand{\Hyp}{{\mathbb H}}

\newcommand{\sone}{{\mathbb S^1}}

\newcommand{\psl}{\PSL(2,\BC)}

\DeclareMathOperator{\B}{B}
\title[Geometry and Rank of Fibered Hyperbolic $3$-Manifolds]{Geometry and Rank of Fibered \\ Hyperbolic $3$-Manifolds}
\author{Ian Biringer}
\begin{document}
\maketitle

\section{Introduction}
 
Recall that the rank of a finitely generated group is the minimal number of elements needed to generate it.  In \cite{white}, M. White proved that the injectivity radius of a closed hyperbolic $3$-manifold $M$ is bounded above by some function of $\rank (\pi_1 (M)) $.  Building on a technique that he introduced, we determine the ranks of the fundamental groups of a large class of hyperbolic $3$-manifolds fibering over the circle.  

Let $\Sigma_g$ be the closed orientable surface of genus $g$ and $\phi : \Sigma_g \to \Sigma_g$ a homeomorphism.  We can construct a $3$-manifold $M_\phi$, the \it mapping torus \rm of $\phi$, as the quotient space $$M_\phi = \Sigma_g \times [0,1] / \sim, \ \ (x,0) \sim (\phi(x),1).  $$  Thurston \cite{thurstonone} has proven that if the map $\phi : \Sigma_g \to \Sigma_g$ is pseudo-anosov then $M_\phi$ can be given a hyperbolic metric.

The fundamental group of $M_\phi$ is given by an HNN-extension $$1 \to \pi_1(\Sigma_g) \to \pi_1(M_\phi) \to \BZ \to 1.$$  Since $\rank (\pi_1 (\Sigma_g)) = 2g  $ it follows that $\rank(\pi_1(M_\phi)) \leq 2g + 1 $.  It is not hard to construct examples where this inequality is strict, but it seems likely that if the gluing map is complicated enough then equality should hold.  As an illustration of this, J. Souto proved in \cite{juan} that given a pseudo-anosov map $\phi : \Sigma_g \to \Sigma_g$, we have for sufficiently large powers $\phi^n$ of $\phi$ that $\rank(\pi_1(M_{\phi^n})) = 2g+1$.  Our main result is the following extension of Souto's theorem.

\begin{theorem}

\label{maintheorem}
Given $\epsilon > 0$ and a closed orientable surface $\Sigma_g$, there are at most finitely many $\epsilon$-thick hyperbolic 3-manifolds $M$ fibering over $S^1$ with fiber $\Sigma_g$ for which $rank( \pi_1(M)) \neq 2g+1$. 
\end{theorem}

Recall that the \it injectivity radius \rm of a hyperbolic manifold $M$, written $\inj (M) $, is defined to be half the length of a shortest homotopically essential loop in $M$, and $M$ is called \it $\epsilon$-thick \rm if $\inj(M) \geq \epsilon$.

\vspace{1mm }

Results similar to Theorem \ref{maintheorem} concerning the Heegaard genus of $M$ are already known; the strongest is due to Bachman and Schleimer, \cite{bs}.  Recall that the \it Heegaard genus \rm of a closed $3$-manifold $M$ is the smallest $g=g (M) $ such that $M$ can be obtained by gluing two genus $g$ handlebodies along their boundaries.  It is easy to see that when $M$ fibers over the circle with fiber $\Sigma_g $ then $g (M) \leq 2g +1 $, and Bachman and Schleimer show that $g (M) =2g +1 $ as long as the monodromy map of $M$ has translation distance at least $2g+1 $ in the curve complex of $\Sigma_g$.  It is likely that the conclusion of Theorem \ref{maintheorem} is true under similar assumptions, but it is not yet clear to us how to prove this.

Before beginning the bulk of this paper, let us sketch the idea behind the proof of Theorem \ref{maintheorem}.  Let $M$ be a hyperbolic $3$-manifold fibering over the circle with fiber $\Sigma_g $.  Following a technique of White \cite{white}, we find a graph $X $ with $\rank (\pi_1 (X)) = \rank (\pi_1 (M)) $ and a $\pi_1 $-surjective mapping $f: X \to M $ whose image has as small length as possible.  We show that if $M$ has large diameter it is most efficient for $X $ to use small edges to fill out the fundamental group of the fiber and a long edge to circumnavigate $M$ in the horizontal direction.  The subgraph of $X $ consisting of all small edges then has rank at least $2g $, since it generates a subgroup of $\pi_1(M) $ isomorphic to $\pi_1 (\Sigma_g) $.  But $X $ must have even larger rank, so $\pi_1(M) = \pi_1(X) \geq 2g +1$.  

The paper is organized as follows.  We begin in Section \ref{prelims} by recalling some standard facts from the theory of Kleinian groups.  In Section \ref{surfacegroups}, we use a lemma of Souto and a compactness argument to provide geometry bounds for certain covers of doubly degenerate hyperbolic manifolds homeomorphic to $\Sigma_g \times \BR $.  The minimal length graphs mentioned above are formally introduced in Section \ref{carriergraphs} and Section \ref{proofoftheorem} contains a proof of Theorem \ref{maintheorem}.  We finish with an appendix that fleshes out a result due to Souto, \cite{rank3}, that gives a convenient decomposition for minimal length $\pi_1 $-surjective graphs in closed hyperbolic $3$-manifolds.

Acknowledgements: I thank Justin Malestein, Nathan Broaddus and Benson Farb for their helpful comments and Juan Souto for many conversations, advice and insight.

\section { Preliminaries }
\label{prelims}

Let $M$ be a hyperbolic $3$-manifold with finitely generated fundamental group.  For the sake of simplicity, we will assume that $M $ has no cusps.  A result of P. Scott in $3$-manifold topology states that $M$ admits a \it compact core, \rm that is a compact submanifold $N $ whose inclusion into $M$ is a homotopy equivalence, \cite{scottcores}.  The connected components of $M \setminus N $ are called the \it ends \rm of $M$.  Marden asked in the 1970s whether $M$ is always homeomorphic to the interior of its compact core; this was recently proven to be true by Agol \cite{agol} and Calegari-Gabai \cite{cg}.  Consequently, if $E $ is an end of $M$ then $E $ is homeomorphic to $\partial E \times [0,\infty) $, where $\partial E $ is the boundary component of $N $ facing $E $.

   Define the \it convex core \rm of $M$ to be the smallest convex submanifold $\CC(M)\subset M $ whose inclusion is a homotopy equivalence.  An end $E $ of $M$ is called \it convex-cocompact \rm if $E \cap \CC(M) $ is compact, and \it degenerate \rm otherwise.  A convex-cocompact end is geometrically a warped product, where the metric on level surfaces of $E \cong \partial E \times [0,\infty) $ grows exponentially with the distance to the boundary of the convex core.  The geometry of degenerate ends is more subtle - we will limit ourselves here to a pertinent example and forward the reader to \cite{thurstonnotes} and \cite{japanese} for the general theory.

\begin{example}
Let $M_\phi$ be the mapping torus of a pseudo-anosov map $\phi: \Sigma_g \to \Sigma_g $.  As mentioned in the introduction, $\pi_1(M_\phi)$ decomposes as $$1 \to \pi_1(\Sigma_g) \to \pi_1(M_\phi) \to \BZ \to 1.$$  Let $N$ be the cyclic cover of $M_\phi$corresponding to the subgroup $\pi_1(\Sigma_g)$.  Then $N$ is homeomorphic to $\Sigma_g \times \BR $, and since it regularly covers a closed manifold we have $\CC(N) =N $, implying that both ends of $N $ are degenerate.  Note that unwrapping a fiber bundle structure for $M_\phi $ gives a product structure $N \cong \Sigma_g \times \BR$ with fibers of bounded diameter, contrasting with the exponential growth of level surfaces in a convex-cocompact end.\end{example}

\subsection{ Simplicial Hyperbolic Surfaces }  We record here some facts about negatively curved surfaces in hyperbolic $3$-manifolds.

 \begin{definition}
Let $M$ be a hyperbolic $3$-manifold.  A \it simplicial hyperbolic surface \rm in $M$ is a map $f: S \to M $, where
\begin{itemize}
\item $S$ is a closed surface equipped with a triangulation $T$

\item $f $ maps each face of $T$ to a totally geodesic triangle in $M$
\item for each vertex $v \in T$ the angles between the $f $-images of the edges adjacent to $v $ sum to at least $2\pi$.
\end{itemize}
\end{definition}

If $f : S \to M$ is a simplicial hyperbolic surface then we get a path-metric on $S$ by requiring that $f$ preserves path lengths.  The metric is smooth and hyperbolic away from the vertices of $T $, at which there are possible excesses of angle.  By the Gauss-Bonnet Theorem, we have $\vol(S) \leq 2\pi | \chi(S) |$.  Bounding the diameter of $S $ by its volume and injectivity radius, we obtain:

\begin{boundeddiameterlemma}[Bonahon]
\label{boundeddiameter}
Assume $f :S \to M$ is an $\epsilon $-thick simplicial hyperbolic surface of genus $g $.  Then $\diam(S) \leq \frac{4}{\epsilon}(2g -2)$.
\end{boundeddiameterlemma}

Mahler's Compactness Theorem (\cite{BP}, E.1) states that the moduli space of $\epsilon $-thick (smooth) hyperbolic surfaces of fixed genus is compact.  Together with the following Proposition, this provides a number of upper bounds on the geometry of $\epsilon $-thick simplicial hyperbolic surfaces, albeit without explicit constants.

 \begin{prop} [Smooth Dominates Simplicial]
\label{domination}
Let $S$ be a closed surface and $d$ a metric on $S $ that is the pullback metric for some simplicial hyperbolic surface.  Then there exists a smooth hyperbolic metric $d_{hyp}$ on $S $ such that for all $x,y \in S $ $$\frac{1 }{ C } d(x,y) \leq d_{hyp} (x,y), $$ where $C >0$ depends only on the topological type of $S$.  Note that if $d$ is $\epsilon $-thick	then $d_{hyp}$ is $\frac{\epsilon}{ C } $-thick.
\end{prop}
\begin{proof}[Proof of Proposition \ref{domination}]
Working in polar coordinates in small neighborhoods around the singular points of $d$, we can explicitly deform $d $ to obtain a smooth metric $d' $ with Gaussian curvature $K \leq -1 $ that is bilipschitz to $d $ with bilipschitz constant depending only on the angles $d $ has around the points in its singular locus.  The argument is very similar to the proof of the $2\pi$-Theorem of Gromov and Thurston \cite{bleiler}, so we will omit it here.  Since the Gauss-Bonnet Theorem gives an upper bound for the sum of these singular angles, $d $ and $d'$ are in fact $C $-bilipschitz for some $C$ depending only on the topological type of $S $.  Define $d_{hyp} $ to be the hyperbolic metric in the conformal class of $d' $.  The Ahlfors-Schwartz Lemma \cite{ahlfors} states that distances measured in $d'$ are less than or equal to distances in $d_{hyp} $; this proves the desired inequality.\end{proof}

As an application, we can use Proposition \ref{domination} and a based version of Mahler's Compactness Theorem to show:

\begin{cor} [Short Markings]
\label{goodmarkings}
Set $\Gamma = \pi_1(\Sigma_g) $ and fix a generating set $X \subset \Gamma $.  Then given $\epsilon,g >0 $ there is a constant $L$ such that whenever $f:S \to M $ is an $\epsilon$-thick simplicial hyperbolic surface of genus $g $ and $p \in S$, there is an isomorphism $\Phi:\Gamma \to \pi_1 (S,p)$ such that the image of each element of $X$ can be represented by a loop based at $p $ of length less than $L $.\end{cor}

\subsection { Algebraic and Geometric Convergence }
Let $\Gamma $ be a finitely generated group and consider a sequence of discrete and faithful representations $\rho_i: \Gamma \to \psl$.  If $(\rho_i) $ converges pointwise to $\rho_\infty: \Gamma \to \psl $, we usually say that $(\rho_i) $ is \it algebraically convergent \rm with $\rho_\infty $ as its algebraic limit.  Alternatively, consider a sequence of subgroups $G_i \subset \psl $; if these converge to a subgroup $G \subset \psl $ in the Hausdorff topology on closed subsets of $\psl $ then we say that $G_i \to G $ \it geometrically. \rm  The case where the two notions of convergence agree is useful enough to warrant additional terminology.  Specifically, if  $\rho_i \to \rho_\infty $ algebraically and $\rho_i (\Gamma) \to \rho_\infty (\Gamma) $ geometrically then one says that $\rho_i \to \rho_\infty $ \it strongly. \rm  

One can interpret the geometric convergence of a sequence of subgroups $G_i \to G_\infty \subset \psl $ in terms of the quotient manifolds $M_i = \Hyp^3 / G_i $.  If we fix a basepoint and baseframe $(p, f) $ for $\Hyp^3 $, for each $i $ we can take the projection $  (p_i, f_i) $ as a basepoint and baseframe for $M_i $.  Then $G_i \to G_\infty $ geometrically if there exist sequences of positive numbers $\epsilon_i \to 0 $ and $R_i \to \infty $, and $(1+\epsilon_i) $-bilipschitz maps $\phi_i : \B(p_i,R_i) \to M_\infty$ sending $(p_i, f_i) $ to $(p_\infty,f_\infty) $. For future reference, we will call the maps $\phi_i $ a sequence of \it almost isometric maps \rm coming from geometric convergence.  Note that using this as our definition, we can speak about a geometrically convergent sequence of framed hyperbolic $3$-manifolds, or even based hyperbolic $3$-manifolds if we forget about the presence of a baseframe.

For a detailed study of algebraic and geometric convergence, see \cite{japanese} and \cite{BP}.

\section{Short Graphs in Doubly Degenerate $\Sigma_g \times \BR $}
\label{surfacegroups}
Assume that $M $ is a hyperbolic $3$-manifold without cusps that is homeomorphic to $\Sigma_g \times \BR $.  Using Waldhausen's Cobordism Theorem \cite{wald}, it is not hard to see that there is an explicit homeomorphism $M \cong \Sigma_g \times \BR $ such that $\CC(M) $ sits inside $M$ as either
\begin{itemize}
\item{ $\Sigma_g \times [0, 1] $, in which case $M$ is convex cocompact}
\item{ $\Sigma_g \times [0, \infty) $, in which case $M$ is called singly degenerate }
\item{ $\Sigma_g \times \BR $, and then $M$ is called doubly degenerate.  }
\end{itemize}  

We mentioned in the introduction that Theorem \ref{maintheorem} is an extension of an earlier theorem of Souto \cite{juan}.  A key ingredient in Souto's proof was the following observation, which is a consequence of the Covering Theorem of Canary and Thurston \cite{canarycovering}.

\begin{lemma}[\cite{juan}] Let $M$ be a doubly degenerate hyperbolic $3$-manifold homeomorphic to $\Sigma_g \times \BR$ and let $\Gamma \subset \pi_1(M) $ be a proper subgroup of rank at most $2g$.  Then $\Gamma$ is free, infinite index and convex-cocompact.
\label{freeconvex}
\end{lemma}

To prove Theorem \ref{maintheorem}, we need an improved version of Lemma \ref{freeconvex} that gives a diameter bound for the convex core of $\Hyp^3 / \Gamma $ in terms of $\inj (M) $ and the length of a set of loops in $M$ generating $\Gamma $.  Our proof will be a compactness argument: we define a topology on the set of wedges of $k$ bounded length loops in $\epsilon $-thick doubly degenerate hyperbolic $3$-manifolds homeomorphic to $\Sigma_g \times \BR $, show that the resulting space is compact and then use continuity to show that there is an upper bound for the corresponding convex core diameters.  

\begin{definition}
Define $\CG = \CG(\epsilon, L, k) $ to be the space of pairs $(M,f) $, where
\begin{enumerate}
\item{ $M$ is a doubly degenerate $\epsilon $-thick hyperbolic $3$-manifold } homeomorphic to $\Sigma_g \times \BR $ 
\item{ $f: \wedge_k \sone \to M $ is an $L $-lipschitz map from the wedge of $k $ circles, endowed with some fixed metric.}
\end{enumerate} 

\noindent We say that $(M_i,f_i) \to (M_\infty,f_\infty) $ if 

\begin{enumerate}
\item{ $(M_i,\star_i) $ converges geometrically to $(M_\infty,\star_\infty) $, where $\star_i $ is the wedge point of $f_i (\wedge_k \sone) $ }
\item{ there is a sequence $\phi_i $ of almost isometric maps coming from the geometric convergence in $(1) $ such that $\phi_i \circ f_i: \wedge_k \sone \to M$ converges pointwise to $f_\infty : \wedge_k \sone \to M_\infty $.  }
\end{enumerate}
\end{definition}

\begin{prop}
\label{compactness}
$\CG$ is compact.

\end{prop}
\begin{proof} 
Let $(M_i,f_i) $ be a sequence in $\CG$ and assume that $\star_i \in M_i $ is the wedge point of $f_i (\wedge_k \sone) $.  For each $i $, Canary's Filling Theorem \cite{canarycovering} gives a simplicial hyperbolic surface in $M_i $ with image passing through $\star_i $.  Using the short markings of these surfaces provided by Corollary \ref{goodmarkings} we can construct representations $\rho_i : \pi_1(\Sigma_g) \to \psl $ with $\Hyp^3 / \rho_i (\Sigma_g) \cong M_i $ so that a fixed base point $\star \in \Hyp^3 $ projects to each $\star_i $ and up to passing to a subsequence, $\rho_i $ converges algebraically to some $\rho_\infty: \pi_1 (\Sigma_g) \to \psl  $.  Since our lower bound on injectivity radius persists through algebraic limits, $\rho_\infty (\pi_1 (\Sigma_g)) $ contains no parabolics.  Work of Thurston and Bonahon\footnote{ Strong convergence follows here from tracing through Thurston's proof of (\cite{thurstonnotes}, 9.2) with the hindsight provided by Bonahon's Tameness Theorem \cite{bonahon}.  A statement of the resulting theorem is given by Canary in (\cite{canarycovering}, 9.1) as a prelude to a series of more general convergence theorems.  } then implies that $\rho_i \to \rho_\infty $ strongly.

Set $M_\infty =\Hyp^3 / \rho_\infty (\Sigma_g)  $ and let $\star_\infty \in M_\infty $ be the projection of $\star $.  Then $(M_i,\star_i) $ converges geometrically to $(M_\infty,\star_\infty) $.  The fundamental group of $M_\infty $ is isomorphic to $\pi_1 (\Sigma_g) $, so Bonahon's Tameness Theorem \cite{bonahon} implies that $M_\infty \cong \Sigma_g \times \BR $.  Moreover, it follows from strong convergence that $M_\infty $ is doubly degenerate.  We can construct a map $f_\infty: \wedge_k \sone \to M_\infty $ by applying Arzela-Ascoli's Theorem to the sequence of maps $\phi_i \circ f_i : \wedge_k \sone \to M_\infty $, where $\phi_i $ is a sequence of almost isometric maps coming from geometric convergence.  Clearly $(M_i,f_i) $ converges to $(M_\infty,f_\infty) $ in $\CG $.\end{proof}

\begin{cor}
\label{boundedconvexcores}
Let $M$ be a doubly degenerate $\epsilon$-thick hyperbolic $3$-manifold homeomorphic to $\Sigma_g \times \BR$ and let $p \in M $ be a basepoint.  Assume that $\Gamma \subset \pi_1(M,p) $ is a proper subgroup that can be generated by $2g $ loops based at $p $ of length less than $L $.  Then $\Gamma $ is convex cocompact and the diameter of the convex core of $\Hyp^3 / \Gamma $ is bounded above by some constant depending only on $L, \epsilon $ and $g$.
\end{cor}

\begin{proof}
Observe that $\Gamma $ determines an element $(M,f) \in \CG =\CG (\epsilon,L,2g) $, with the extra property that $f $ is not $\pi_1 $-surjective.  The subset of $\CG $ consisting of pairs $(M,f) $ for which $f $ is not $\pi_1 $-surjective is closed in $\CG$, and therefore compact by Proposition \ref{compactness}.  Lemma \ref{freeconvex} implies that for all such $(M,f) $ the cover $M_{\pi_1 f} $ of $M$ corresponding to the $\pi_1 $-image of $f $ is convex cocompact.  It is not hard to see that if $(M_i,f_i) \to (M_\infty, f_\infty) \in \CG $ then $(M_i)_{\pi_1 f_i} \to (M_\infty)_{\pi_1 f_\infty}$ algebraically after picking appropriate markings.  The diameter of the convex core of a hyperbolic $3$-manifold is continuous with respect to algebraic convergence, so the diameter of the convex core of $(M_i)_{\pi_1 f_i} $ varies continuously over $\CG$.  This proves the claim.\end{proof}

\section{Carrier Graphs}
\label{carriergraphs}

In the following, assume $M$ is a closed hyperbolic $3$-manifold.

\begin{definition}  A \it carrier graph \rm for $M$ is a graph $X$ and a map $f: X \to M$ which induces a surjection on fundamental groups.
\end{definition}

\noindent \it Standing Assumption: \rm In this paper we are interested in generating sets of minimal size, which correspond to carrier graphs with $\rank(\pi_1(X))=\rank(\pi_1(M))$.  From now on all carrier graphs will be assumed to have this property.
\vspace{2mm}

If a carrier graph $f:X \to M$ is rectifiable, we can pull back path lengths in $M$ to obtain an pseudo-metric on $X$.  Collapsing to a point each zero-length segment in $X$ yields a new carrier graph with an actual metric; from now on we will assume all carrier graphs are similarly endowed.  Define the \it length \rm of a carrier graph to be the sum of the lengths of its edges, and a \it minimal length carrier graph \rm to be a carrier graph which has smallest length (over all carrier graphs of minimal rank).  An argument using Arzela-Ascoli's Theorem, \cite{white}, shows that minimal length carrier graphs exist in any closed hyperbolic $3$-manifold.

The following Proposition shows that minimal length carrier graphs are geometrically well behaved.

\begin{prop}[White, \cite{white}]
\label{minimallengths}
Assume $f : X \to M$ is a minimal length carrier graph in a closed hyperbolic $3$-manifold $M$.  Then $X$ is trivalent with $2(\rank(\pi_1(M)) - 1)$ vertices and $3 (\rank(\pi_1(M)) - 1)$ edges, each edge in $X$ maps to a geodesic segment in $M$, the angle between any two adjacent edges is $\frac{2\pi}{3}$, and the image of any simple closed path in $X$ is an essential loop in $M$.
\end{prop}

We conclude this section with a technical result that is instrumental in our proof of Theorem \ref{maintheorem}.  A slightly more general theorem was proven by Souto in \cite{rank3}, but the proof given there is somewhat incomplete.   We include a full proof of the more general result in Appendix \ref{proofofchains}.

\begin{prop}[Chains of Bounded Length]
\label{chainsofboundedlength}
Let $M$ be a closed hyperbolic 3-manifold with $f: X \to M$ a minimal length carrier graph. Then we have a sequence of (possibly disconnected) subgraphs $$\emptyset = Y_0 \subset Y_1 \subset \ldots \subset Y_k = X$$ such that the length of any edge in $Y_{i+1} \setminus Y_i$ is bounded above by some constant depending only on $\inj(M)$, $\rank(\pi_1(M))$, $\length(Y_i)$ and the diameters of the convex cores of the covers of $M$ corresponding to $f_*(\pi_1(Y_i^j))$, where $Y_i^1, \ldots, Y_i^n$ are the connected components of $Y_i$.
\end{prop}

\section{Proof of Theorem \ref{maintheorem}}
\label{proofoftheorem}

Fix $\epsilon,g>0$ and assume that $M $ is an $\epsilon$-thick hyperbolic $3$-manifold fibering over the circle with fiber $\Sigma_g$.  The goal of this section is to prove that there are only finitely many such $M$ for which $\rank (\pi_1(M)) \neq 2g +1 $.  We begin, however, with a quick computation concerning $M$'s girth.

\begin{definition}
\label{circumference}
The \it waist length \rm of $M$, denoted $\waist (M) $, is the smallest length of a loop in $M$ that projects to a generator of $\pi_1(S^1)$.  
\end{definition}

\begin{prop} [Fibered $3$-Manifolds Have High BMI]
\label{waistversusdiameter}
Let $M$ be an $\epsilon$-thick hyperbolic $3$-manifold fibering over the circle with fiber $\Sigma_g$.  Then $$2 \diam(M) -\frac{16}{\epsilon}(2g -2) \leq \waist(M) \leq 2 \diam(M).  $$
\end{prop}

\begin{proof}

Assume that $\gamma $ is a loop realizing the waist length of $M$.  Canary's Filling Theorem \cite{canarycovering} implies that every point in the cyclic cover of $M$ corresponding to the fundamental group of the fiber lies in the image of a simplicial hyperbolic surfaces for which the inclusion map is a homotopy equivalence.  Projecting down, this provides an exhaustion of $M$ by simplicial hyperbolic surfaces in the homotopy class of the fiber.  By homological considerations, any such surface must intersect $\gamma $.  The Bounded Diameter Lemma (see Section \ref{prelims}) then implies that $\diam(M) \leq \frac{1}{2}\waist (M) +\frac{8}{\epsilon}(2-2g)$.  This establishes the first inequality.

For the second, recall that the fundamental group of $M$ is generated by the set of all loops in $M$ of length less than $2 \diam(M).  $  Any generating set for $\pi_1(M) $ must contain a loop that encircles $M$'s waist, so the waist length of $M$ is at most twice its diameter.  \end{proof}

There are only finitely many hyperbolic $3$-manifolds with diameter less than a given constant.  Proposition \ref{waistversusdiameter} then gives a similar finiteness result for thick hyperbolic $3$-manifolds fibering over the circle with a fixed fiber and bounded waist length.

\vspace{2mm }

We are now ready to prove the main result of this note.

\begin{maintheoremenv}
Given $\epsilon,g > 0$ there are at most finitely many $\epsilon$-thick hyperbolic 3-manifolds $M$ fibering over $\sone$ with fiber $\Sigma_g$ for which $rank( \pi_1(M)) \neq 2g+1$.
\end{maintheoremenv}

\begin{proof}
Assume that $M$ is an $\epsilon$-thick hyperbolic $3$-manifold fibering over the circle with fiber $\Sigma_g$ and $\rank( \pi_1(M)) \leq 2g$.  We will show that the waist length of $M$ is bounded by some constant depending only on $\epsilon$ and $g$. 

Let $f : X \to M$ be a minimal length carrier graph.  By Proposition \ref{chainsofboundedlength}, there is a constant $L$ and a chain of (possibly disconnected) subgraphs $$\emptyset = Y_0 \subset Y_1 \subset \ldots \subset Y_k = X$$ with $\length (Y_{i +1}) $ bounded above by some constant depending only on $\epsilon $, $g$, $\length(Y_i)$ and the diameters of the convex cores of the covers of $M$ corresponding to the fundamental groups of the connected components of $Y_i $.  

Assume for the moment that no connected component of $Y_i $ runs all the way around $M$'s waist, so that each lifts homeomorphically to the cyclic cover $M_{ \pi_1 (\Sigma_g) } $ of $ M$.  Since $\rank( \pi_1(X)) \leq 2g$, the components of $Y_i $ have even smaller rank and thus cannot generate the fundamental group of $M_{ \pi_1 (\Sigma_g) } $.  Therefore Corollary \ref{boundedconvexcores} applies to bound the diameters of the associated convex cores in terms of $\length(Y_i)$, $\epsilon $ and $g $.  It follows that $\length (Y_{i +1}) $ is also bounded above by $\length(Y_i)$, $\epsilon $ and $g $.  

Applying this argument iteratively, we obtain a length bound for the first subgraph $Y_i $ that has a component which navigates the waistline of $M$.  The length bound depends on $\epsilon $, $g $ and the index of the subgraph, but since there are at most $3(\rank (\pi_1(M)) - 1) $ edges in $X $ the number of subgraphs in our chain is also limited.  Therefore we have that the waist length of $M$ is bounded by a function of $\epsilon $ and $g $.\end{proof}

Under slight modifications, the proof of Theorem \ref{maintheorem} shows that for mapping tori with large waist length there is only one Nielsen equivalence class of minimal size generating sets for $\pi_1(M)$.  The interested reader may compare our proof with \cite{rank3} for more details.

% Second, the number of exceptional cases in Theorem \ref{maintheorem} is in principle computable, but the proof above does not illustrate this.  However, the only step of our proof that does not yield explicit constants is the compactness argument in Section \ref{surfacegroups}.  Recall that at the onset of this argument we had an $\epsilon $-thick doubly degenerate hyperbolic $3$-manifold $M \cong \Sigma_g \times \BR $ and a proper subgroup $\Gamma \subset \pi_1(M) $ generated by $2g $ loops in $M$ of bounded length.  Our goal was to provide a diameter bound for the convex core of $M_\Gamma = \Hyp^3/ \Gamma $.  This can be accomplished more explicitly as follows.  First, recall from the proof of Lemma \ref{freeconvex} that $M_\Gamma $ is a handlebody.  Then $\partial \CC (M_\Gamma) $ is $\epsilon $-thick, for the lower bound on injectivity radius rules out short curves that are homotopically essential in $M_\Gamma $ and the length bound on the generators prevents short compressible curves.  Whose purpose was to give diameter bounds for the convex cores of covers of $\epsilon $-thick doubly degenerate hyperbolic manifolds $M \cong Sigma_g \times \BR $ whose fundamental groups are $ generated by $2g $ short loops in $M$.  This can be accomplished more explicitly as follows.

\appendix
\section{Chains of Bounded Length }

\label{proofofchains}

We prove here the generalization of Proposition \ref{chainsofboundedlength} promised in Section \ref{carriergraphs}.  The idea of the proof given below was originally sketched by Souto in \cite{rank3}; the purpose of this Appendix is to fill in some missing details.

Assume that $M =\Hyp^3 / \Gamma $ is a closed hyperbolic $3$-manifold and $f: X \to M $ is a minimal length carrier graph. Choose an edge $e \subset X $ and a subgraph $Y \subset X $.  Our first goal will be to provide a useful definition of the length of $e $ relative to the subgraph $Y $.  This should vanish when $e \subset Y $ and should agree with the hyperbolic length of $f (e) $ when neither of the vertices of $e$ lies inside $Y $.  If $X $ is embedded as a subset of $M$ with $f $ the inclusion map, then relative length is similar to the length $e$ has outside of the hyperbolic convex hulls of the components of $Y $ that $e $ touches, but we need to do our measurements in the universal cover and throw out sections of $e$ that lie inside some of the thin parts of $M$.  

To clarify this, fix a universal covering $\pi_X:\tilde { X }\to X $ and a lift $\tilde { f }:\tilde { X } \to \Hyp^3 $ of $f $.  Assume that a vertex $v$ of $e $ lies in a connected component $Z_v \subset Y $ and choose lifts $\hat { e },\hat { Z_v } \subset \tilde { X } $ of $e$ and $Z $ that touch above $v$.  Let $\Gamma_{\tilde { f } (\hat { Z_v }) } $ be the subgroup of $\Gamma  $ that leaves $\tilde { f } (\hat { Z_v }) $ invariant.

\begin{definition}[Thick Convex Hulls]

The \it thick convex hull of $\tilde { f } (\hat { Z_v })  $, \rm written $\TCH(\tilde { f } (\hat { Z_v })  )$, is the smallest convex set $K $ containing $\tilde { f } (\hat { Z_v })   $ such that for every $\gamma \in \Gamma_{\tilde { f } (\hat { Z_v })  } $ and $x \in \Hyp^3 \setminus K$, we have $\dist(\gamma(x),x) \geq 1$. 
\end{definition}

\begin{definition}[Edge Length Relative to a Subgraph]
Define the \it length of $e$ relative to $Y$, \rm denoted $\length_Y (e) $, to be the length of the part of $\tilde{f} (\hat {e }) $ that lies outside of $\TCH(\tilde { f } (\hat { Z_v })) $ for each vertex $v $ of $e $ contained in $Y $.  

\end{definition}

It is easy to see that the relative length of $ e $ is well-defined, independent of the lifts chosen above.  The definition is a bit less complicated if we assume that $X $ is embedded as a subset of $M$.  For then we can lift $e $ directly to $\Hyp^3 $ along with any connected components of $Y $ that $e $ touches, and then measure the length of $e$'s lift outside of the thick convex hulls of the lifted subgraphs.  In the proofs below, we will assume $X $ to be embedded in order to remove a level of notational hinderance.  The arguments will be exactly the same in the general case.

Although an edge can have very long absolute length while having short length relative to a subgraph $Y $, we can bound this difference if we have some control over the geometry of the covers of $M$ corresponding to the fundamental groups of the components of $Y $.

\begin{lemma}
\label{relativeversusabsolute}
Assume that $M$ is a closed hyperbolic $3$-manifold, $f: X \to M $ is a minimal length carrier graph, $Y $ is a subgraph of $X $ and $e $ is an edge of $X \setminus Y $.  Then $\length (e) $ is bounded above by a constant depending only on $\length_{ Y } (e) $, $\length (Y) $, $\inj (M) $, $\rank (\pi_1 (M)) $ and the diameters of the convex cores of the covers of $M$ corresponding to the components of $Y $ that $e $ touches.
\end{lemma}

\begin{proof}
As mentioned above, we forget about $f $ and assume that $X $ is embedded as a subset of $M$.  Suppose that $e $ shares a vertex with a connected component $Z \subset Y $, and let $\tilde { e },\tilde { Z } \subset \Hyp^3 $ be lifts that touch above	that vertex.  It suffices to show that the Hausdorff distance from $\tilde { Z } $ to $\TCH (\tilde { Z }) $ is bounded by the quantities mentioned in the statement of the Lemma.  For since $X $ is minimal length, $\tilde { e } \cap \TCH(\tilde { Z }) $ must minimize the distance from $\tilde { e } \cap \partial \TCH (\tilde { Z }) $ to $\tilde { Z } $; thus a bound on the Hausdorff distance between $\tilde { Z } $ and $\TCH (\tilde { Z }) $ limits the length that $\tilde{e } $ can have inside of $\TCH (\tilde { Z }) $.

We first claim that the hyperbolic distance from $\tilde { Z } $ to $\CH (\Lambda (\Gamma_{\tilde { Z } })) $ is bounded above by a constant depending only on $\inj (M) $ and $\rank (\pi_1(M)) $.  Choose an infinite piecewise geodesic path $\gamma \subset \tilde { Z } $ that projects to a simple closed curve in $Y $ and let $g \in \Gamma_{\tilde { Z } } $ be the corresponding deck transformation.  Taking a maximal sequence of consecutive edges of $\gamma $ that project to distinct edges in $M $ yields a subpath $\gamma' $ whose $g$-translates cover $\gamma $.  Note that the orthogonal projection of $\gamma' $ to $\axis (g) $ has length equal to the translation distance of $g$, which is at least $\inj (M) $.  By Lemma \ref{minimallengths}, $X $ has $3 (\rank(\pi_1(M)) - 1) $ edges; the number of edges in $\gamma' $ can certainly be no greater than this.  Thus there is an edge of $\gamma $ whose orthogonal projection to $\axis (g) $ has length at least $\frac{\inj (M) }{ 3 (\rank(\pi_1(M)) - 1) } $.  It follows from elementary hyperbolic geometry that there is a point on this edge whose distance from $\axis (g) $ is bounded above by a constant depending on that length; this proves the claim.

Now $\tilde { Z } $ and $\CH (\Lambda (\Gamma_{\tilde { Z } })) $ are both invariant under the action of $\Gamma_{\tilde { Z } } $ with quotients of bounded diameter, so our limit on the hyperbolic distance between them translates into a bound on their Hausdorff distance.  But if $\tilde { Z } $ is Hausdorff-close to a convex set then it must also be Hausdorff-close to its convex hull, $\CH (\tilde { Z }) $.  Since the Hausdorff distance from $\CH (\tilde { Z }) $ to $\TCH (\tilde { Z }) $ is controlled by $\inj (M) $, we have a bound on the Hausdorff distance between $\tilde { Z } $ and $\TCH (\tilde { Z }) $ .
\end{proof}

\vspace{2mm}

For a subgraph $Z \subset X$, we define the \it length of $Z$ relative to $Y$ \rm to be $$\length_Y(Z) = \sum_{\text{edges $e \subset Z$}} \length_Y(e).  $$

Using our definition of relative length, we can streamline the formulation of Proposition \ref{chainsofboundedlength}.  The statement given earlier follows from this one after applying Lemma \ref{relativeversusabsolute}.

\begin{prop}[Chains of Bounded Length]
There is a universal constant $L $ with the property that if $M$ is a closed hyperbolic 3-manifold and $f: X \to M$ is a minimal length carrier graph then we have a sequence of (possibly disconnected) subgraphs $$\emptyset = Y_0 \subset Y_1 \subset \ldots \subset Y_k = X$$ such that $\length_{Y_i}(Y_{i+1}) < L$ for all $i$.\end{prop}

\begin{proof}

It is a standard fact in hyperbolic geometry that there exist a universal constant $C >0 $ with the following property: 

\begin{enumerate}
\item{ any path in $\Hyp^3 $ made of geodesic segments of length at least $C$ connected with angles at least $\frac{\pi}{3}$ is a quasi-geodesic.  }
\end{enumerate}

\noindent There is also a constant $D > C$ such that

\begin{enumerate}
\item[(2)]{if $N \subset \Hyp^3 $ contains the axis of a hyperbolic isometry $\gamma $ and $\dist(x,\gamma(x)) \geq 1 $ for all $x \in \Hyp^3 \setminus N $, then $\dist(x,\gamma (x)) \geq C$ for all $x \in \Hyp^3 \setminus \CN_D (N) $, } 
\item[(3)]{any geodesic ray emanating from a convex subset $K \subset \Hyp^3 $ that leaves $\CN_D (K)$ meets $\partial \CN_D(K) $ in an angle of at least $\frac{\pi}{3}$, }

\end{enumerate}

\noindent and finally a constant $B >0 $ for which

\begin{enumerate}
\item[(4)]{ any geodesic exiting a convex subset $K \subset \Hyp^3 $ will exit $\CN_D (K)$ after an additional length less than $B $.  }
\end{enumerate}

We will show that if $Y $ is any subgraph of $X $ then there is an edge in $X \setminus Y$ of length at most $L = C + 2B$ relative to $Y$; applying this iteratively will give the chain of subgraphs in the statement of the Proposition.

So, suppose that $Y $ is a subgraph of $X $.  Observe that since the fundamental group of a closed hyperbolic manifold cannot be free, there is an essential closed loop $\gamma \subset X $ that is nullhomotopic in $M$.  Furthermore, since $\pi_1(M) $ does not split as a free product, \cite{hempel}, we can pick $\gamma $ so that it has no subpath contained entirely in $Y$ that is also a closed loop nullhomotopic in $M$.  Lifting $\gamma $ to $\Hyp^3 $ then gives a closed loop $\tilde{\gamma} \subset \Hyp^3 $ such that each time $\tilde{\gamma} $ touches a component of $\pi_M^{ - 1} (Y) $ it enters and leaves that component using different edges of $\pi_M^{ - 1} (X \setminus Y) $.

Consider a maximal segment of $\tilde{\gamma}$ that is contained in a component ${\tilde { Z } } $ of $ \pi_M^{ - 1} (Y) $ and let $e$ and $f $ be the edges that $\tilde{\gamma}$ traverses before and after the segment in ${\tilde { Z } } $.  If $e$ or $f$ has length less than $L $ relative to $Y $, then we are done.  Otherwise, the two edges have a length of at least $L $ left after exiting $ \TCH({\tilde { Z } }) $, so by $(4)$ both of these edges must exit $\CN_D (\TCH({\tilde { Z } })) $; let $e_0$ and $f_0$ be the points where they meet $\partial \CN_D (\TCH({\tilde { Z } }))  $.  Assume for the moment that the distance between $e_0 $ and $f_0 $ is less than $C $.  Then by $(2) $, $e$ and $f$ project to different edges in $X $.  Substituting $\pi_M (e \cap \CN_D (\TCH({\tilde { Z } }))) \subset X $ with the projection of the geodesic between $e_0 $ and $f_0 $ therefore yields a new carrier graph for $M$, and since the new edge has length less than $C $ while the old has length at least $D $ our new carrier graph has shorter length than $X $.  This contradicts the minimality of $X $, so $\dist (e_0,f_0) \geq C $.

We can now create a new closed path in $\Hyp^3 $ from $\tilde{\gamma } $ as follows: each time $\tilde{\gamma } $ traverses a component ${\tilde { Z } } $ of $ \pi_M^{ - 1} (Y) $, replace the part of $\tilde {\gamma } $ that lies inside $\CN_D (\TCH({\tilde { Z } })) $ by the geodesic with the same endpoints.  Then the new path is composed of geodesic segments of length at least $C$, and by $(3)$, the segments intersect with angles at least $\frac{\pi}{3} $.  Therefore it is a quasi-geodesic.  Since it is also closed, this is impossible.  
\end{proof}

 \end{document}